\input amstex\documentstyle{amsppt}  
\pagewidth{12.5cm}\pageheight{19cm}\magnification\magstep1
\topmatter
\title Remarks on totally positive flag manifolds\endtitle
\author G. Lusztig\endauthor
\address{Department of Mathematics, M.I.T., Cambridge, MA 02139}\endaddress
\thanks{Supported by NSF grant DMS-2153741}\endthanks
\endtopmatter   
\document

\define\si{\sim}

\define\sqc{\sqcup}

\define\qua{\quad}

\define\part{\partial}

\define\m{\mapsto}

\define\sub{\subset}    

\define\T{\times}
\define\ti{\tilde}
\define\nl{\newline}
\redefine\i{^{-1}}

\define\a{\alpha}
\redefine\b{\beta}
\redefine\c{\chi}
\define\g{\gamma}

\redefine\o{\omega}

\define\s{\sigma}
\redefine\t{\tau}

\redefine\l{\lambda}
\define\z{\zeta}
\define\x{\xi}

\define\Th{\Theta}

\define\CC{\bold C}

\define\NN{\bold N}

\define\PP{\bold P}

\define\RR{\bold R}

\define\cb{\Cal B}

\define\cp{\Cal P}

\define\ct{\Cal T}

\define\tR{\ti R}
\define\tS{\ti S}

\head Introduction\endhead
\subhead 0.1\endsubhead
We often identify an algebraic variety defined over $\RR$ with its set
of $\RR$-points. Let $G$ be a split, connected, simply connected semisimple algebraic
group of simply laced type defined over $\RR$ with a fixed pinning
$(B^+,B^-,x_i,y_i (i\in I))$ as in \cite{L94}. Here $B^+,B^-$ are
opposed Borel subgroups of $G$ defined over $\RR$ with unipotent
radicals $U^+,U^-$ and $x_i:\RR@>>>U^+,y_i:\RR@>>>U^-$ are certain
imbeddings; $I$ is a finite set. Let $T=B^+\cap B^-$. Let
$G_{>0},U_{>0}^+,U_{>0}^-,T_{>0}$ be the (open) sub-semigroups of
$G,U^+,U^-,T$ defined in \cite{L94}.

For any $(u,u')\in U^-_{>0}\T U^-_{>0}$ we have defined in
\cite{L94, 7.1} the open subset
$$\ct_{u,u'}=\{t\in T_{>0};tut\i u'{}\i\in U^-_{>0}\}$$
of $T_{>0}$ and proved that it is nonempty. In this paper we will
state a conjecture on the structure of $\ct_{u,u'}$ (see 0.3), we
prove some special cases of it and we derive some consequences of it.

\subhead 0.2\endsubhead
We introduce some notation. 
For any $\l=(\l_i)_{i\in I}\in\NN^I$ let $V_\l$ be an irreducible rational 
representation of $G$ (over $\CC$) whose highest weight is $\l$. We fix a highest weight 
vector $e_\l$ of $V_\l$. 
Let $\b_\l$ be the canonical basis of $V_\l$ that contains $e_\l$, see \cite{L90}.
Let $e'_\l$ be the lowest weight vector in $\b_\l$. 

Let $j\in I$ and let $\l=\o(j)\in\NN^I$ be such that $\l_i=1$ if $i=j$, $\l_i=0$
if $i\ne j$ (a fundamental weight). If $u\in U^-$ we can
write $ue_{\o(j)}$ as an $\RR$-linear combination of vectors in $\b_{\o(j)}$; let
$Z_j(u)\in\RR$ be the coefficient of $e'_{\o(j)}$ in this linear combination.
This defines a function $Z_j:U^-@>>>\RR$ (in fact a morphism of real
algebraic varieties). From \cite{L90}, \cite{L94} it is known that
$Z_j(U_{>0}^-)\sub\RR_{>0}$. Hence
for $(u,u',t)\in U_{>0}^-\T U_{>0}^-\T T_{>0}$ such that
$t\in\ct_{u,u'}$,
$$z(u,u',t)=(Z_j(tut\i u'{}\i))_{j\in I}\in\RR_{>0}^I$$
is defined.

We define
$$\Th:\{(u,u',t)\in U_{>0}^-\T U_{>0}^-\T T_{>0};t\in\ct_{u,u'}\}
@>>>U_{>0}^-\T U_{>0}^-\T\RR_{>0}^I$$
by $(u,u',t)\m(u,u',z(u,u',t))$.

We can now state:

\proclaim{Conjecture 0.3} $\Th$ is a homeomorphism.
\endproclaim
This will be verified in \S1 in some examples.

\subhead 0.4\endsubhead
Let $\cb$ be the (real) flag manifold of $G$ that is, the real algebraic
variety whose points are the Borel subgroups of $G$ defined over $\RR$.
Following \cite{L4,\S8}, let $\cb_{>0}$ be the (open) subset
$$\{u'B^+u'{}\i;u'\in U_{>0}^-\}=\{uB^-u\i;u\in U_{>0}^+\}$$
of $\cb$ (the last equality is proved in \cite{L94, 8.7}). According to
\cite{L94, 8.9},

(a) {\it for any $g\in G_{>0}$ there is a unique $B\in G_{>0}$ such that
$g\in B$;}
\nl
moreover, the map $\z:G_{>0}@>>>\cb_{>0}$ given by
$g\m B$ is continuous. According to \cite{L21, 5.5(a)}, the map $\z$ is
surjective. Thus, the fibres of $\z$ (that is, the sets $B\cap G_{>0}$
for various $B\in\cb_{>0}$) define a partition of $G_{>0}$ into non-empty closed
subsets indexed by $\cb_{>0}$. 

We are interested in the study of the open set
$B\cap G_{>0}$ of $B$ (for any $B\in\cb_{>0}$). The following result
(conjectured in \cite{L21,\S5}) shows (assuming 0.3) that this open set is
homeomorphic to a product of copies of $\RR_{>0}$.

\proclaim{Proposition 0.5} Assume that 0.3 holds for $G$. For any
$B\in\cb_{>0}$ there exists a canonical homeomorphism
$$\s_B:B\cap G_{>0}@>\si>>U_{>0}^+\T\RR_{>0}^I.$$
\endproclaim
The proof is given in \S2.

\subhead 0.6\endsubhead
We now fix $J\sub I$. Let $P_J^+$ be the subgroup of $G$ generated by $B^+$
and by $\{y_j(a);j\in J,a\in\RR\}$. Let $\cp^J$ be the set of
subgroups of $G$ that are $G$-conjugate to $P_J^+$. Following \cite{L98} we 
define $\cp^J_{>0}$ to be the set of subgroups $P\in\cp^J$ such that
$\g_P:=\{B\in\cb_{>0};B\sub P\}$ is nonempty. 
The following result is a generalization of 0.4(a).

\proclaim{Proposition 0.7}Let $g\in G_{>0}$. There is a unique $P\in\cp^J_{>0}$
such that $g\in P$.
\endproclaim
The proof is given in 3.3.

\subhead 0.8\endsubhead
From 0.7 we see that there is a well defined map $\z_J:G_{>0}@>>>\cp^J_{>0}$ given by
$g\m P$ where $P\in\cp^J_{>0}$ contains $g$. It generalizes the map $\z:G_{>0}@>>>\cb_{>0}$
in 0.4. It is again continuous. It is also surjective (this follows from the surjectivity of
$\z$).  Thus, the fibres of $\z_J$ (that is, the sets $P\cap G_{>0}$
for various $P\in\cp^J_{>0}$) define a partition of $G_{>0}$ into non-empty closed
subsets indexed by $\cp^J_{>0}$. 
Note that, if $P\in\cp^J_{>0}$, then $P\cap G_{>0}=\sqc_{B\in\g_P}(B\cap G_{>0})$.
(In other words, if $g\in P\cap G_{>0}$ and if $B\in\cb_{>0}$ is defined by
$g\in B$ then $B\sub P$. Indeed, we have $B\sub P'$ for a unique $P'\in\cp^J_{>0}$ so that
$g\in P',g\in P$ and then $P=P'$ by 0.7.)
Using then 0.5 (under the assumption that 0.3 holds) we see that there is a well
defined bijection

(a) $P\cap G_{>0}=\g_P\T U_{>0}^+\T\RR_{>0}^I$
\nl
whose restriction to $B\cap G_{>0}$ (for any $B\in\g_P$) is given by
$g\m(B,\s_B(g))$. From the definitions we see that 

(b) {\it the bijection (a) is a homeomorphism.}
\nl
In 3.4 we show that

(c) {\it if $P\in\cp^J_{>0}$ then $\g_P$ is homeomorphic to a product of copies of 
$\RR_{>0}$.}
\nl
Combining (b),(c) we obtain:

\proclaim{Proposition 0.9}  Assume that 0.3 holds for $G$. For any
$P\in\cp^J_{>0}$, the intersection $P\cap G_{>0}$ is homeomorphic to
a product of copies of $\RR_{>0}$.
\endproclaim

\head 1. Examples\endhead
\subhead 1.1\endsubhead
In this section we shall give some examples when 0.3 holds.
We first assume that $G=SL_2(\RR)$ with the standard pinning. Let
$$u=\left(\matrix 1&0\\a&1\endmatrix\right)\in U_{>0}^-,\qua
u'=\left(\matrix 1&0\\a'&1\endmatrix\right)\in U_{>0}^-,$$
so that $(a,a')\in\RR_{>0}^2$.
Now $T_{>0}$ is the set of all matrices
$$\left(\matrix r&0\\0&s\endmatrix\right)$$
with $(r,s)\in\RR_{>0},rs=1$; we can identify $T_{>0}$ with
$\RR_{>0}$ by $(r,s)\m R=s/r$.
Let $t_R$ be the element of $T_{>0}$ corresponding to
$R\in\RR_{>0}$. If $t=t_R$, then
$$tut\i u'{}\i=\left(\matrix1&0\\Ra-a'&1\endmatrix\right)$$
so that
$$\ct_{u,u'}=\{R\in\RR; Ra-a'>0\}.$$
Now $\t:\ct_{u,u'}@>>>\RR_{>0}$, $R\m Ra-a'$ is a homeomorphism
$\ct_{u,u'}@>>>\RR_{>0}$.
This shows that 0.3 holds in our case: the map $Z_i:U^-@>>>\RR$ (for
the unique $i\in I$) attaches to a matrix 
$$\left(\matrix 1&0\\x&1\endmatrix\right)\in U^-$$
the number $x$.

\subhead 1.2\endsubhead
In this subsection we assume that $G=SL_3(\RR)$ with the standard
pinning. We shall prove 0.3 in this case. Let
$$u=\left(\matrix 1&0&0\\a&1&0\\c&b&1\endmatrix\right)\in U_{>0}^-,\qua
u'=\left(\matrix 1&0&0\\a'&1&0\\c'&b'&1\endmatrix\right)\in U_{>0}^-,$$
so that $(a,b,c,ab-c,a',b',c',a'b'-c')\in\RR_{>0}^8$.
Now $T_{>0}$ is the set of all matrices
$$\left(\matrix r&0&0\\0&s&0\\0&0&p\endmatrix\right)$$
with $(r,s,p)\in\RR_{>0},rsp=1$; we can identify $T_{>0}$ with
$\RR_{>0}^2$ by $(r,s,p)\m(R,S)=(s/r,p/s)$.
Let $t_{R,S}$ be the element of $T_{>0}$ corresponding to
$(R,S)\in\RR_{>0}^2$. If $t=t_{R,S}$, then
$$tut\i u'{}\i=\left(\matrix1&0&0\\Ra-a'&1&0\\RSc-Sa'b+a'b'-c'&Sb-b'&1
\endmatrix\right)$$
so that
$$\align&\ct_{u,u'}=\{(R,S)\in\RR^2; Ra-a'>0,Sb-b'>0,\\&
(Rc-a'b)S+a'b'-c'>0, R(S(ab-c)-ab')+c'>0\}.\tag a\endalign$$
(The last inequality is obtain by rewriting
$(Ra-a')(Sb-b')-RSc+Sa'b-a'b'+c'>0$.) Note that any $(R,S)$ in the right
hand side of (a) automatically satisfies $R>0,S>0$ since $Ra>a',Sb>b'$.
We define $\t:\RR^2@>>>\RR^2$ by
$$\t(R,S)=((Rc-a'b)S+a'b'-c',R(S(ab-c)-ab')+c').$$
Let $(A,B)\in\RR_{>0}^2$. We show that $\t\i(A,B)$ consists of two
elements. Let $(R,S)\in\t\i(A,B)$ that is 
$$RSc=a'bS-a'b'+c'+A,RS(ab-c)=ab'R-c'+B.$$
We have
$$(ab-c)(a'bS-a'b'+c'+A)=c(ab'R-c'+B)$$
so that
$$R=(ab-c)a'b(ab'c)\i S+(-aa'bb'+a'b'c+abc'+(ab-c)A-cB)(ab'c)\i.\tag b$$
Substituting this into $RSc=a'bS-a'b'+c'+A$ we obtain
$$(ab-c)a'bS^2+(-2aa'bb'+a'b'c+abc'+(ab-c)A-cB)S+(a'b'-c'-A)ab'=0.$$
We set $\tR=R-a'/a,\tS=S-b'/b$. We obtain
$$\align&(ab-c)a'b(\tS+b'/b)^2+(-2aa'bb'+a'b'c+abc'+(ab-c)A-cB)(\tS+b'/b)\\&
+(a'b'-c'-A)ab'=0\endalign$$
that is
$$(ab-c)a'b\tS^2+(abc'-a'b'c+(ab-c)A-cB)\tS-(A+B)cb'/b=0;\tag c$$
Since $(-(A+B)cb'/b)/((ab-c)a'b)\i<0$, the equation (c) for $\tS$ has two
roots, $\tS_+,\tS_-$ such that $\tS_+\in\RR_{>0},\tS_-\in-\RR_{>0}$.
Thus we have 
$$\tS\in\{\tS_+,\tS_-\}$$
and
$$\tS_++\tS_-=-\mu((ab-c)a'b)\i,\tag d$$
where $\mu=abc'-a'b'c+(ab-c)A-cB$.

Note that
$$\tS_+=(-\mu+\sqrt{\mu^2+4(ab-c)a'b'c(A+B)})/(2(ab-c)a'b),$$
$$\tS_-=(-\mu-\sqrt{\mu^2+4(ab-c)a'b'c(A+B)})/(2(ab-c)a'b),$$
where $\sqrt:\RR_{>0}@>>>\RR_{>0}$ is the square root.

In the case where $\tS=\tS_+$ (resp. $\tS=\tS_-$) we set $\tR=\tR_+$
(resp. $\tR=\tR_-$). We can rewrite (b) as
$$\tR_{\pm}=(ab-c)a'b(ab'c)\i\tS_{\pm}\mu(ab'c)\i.$$
Using (d) this becomes
$$\tR_{\pm}=(ab-c)a'b(ab'c)\i\tS_{\pm}-((ab-c)a'b)(ab'c)\i(\tS_++\tS_-)$$
that is
$$\tR_{\pm}=(ab-c)a'b(ab'c)\i(\tS_{\pm}-\tS_+-\tS_-).$$
Thus, 
$$\tR_+=-(ab-c)a'b(ab'c)\i\tS_-\in\RR_{>0}$$
(since $\tS_-\in-\RR_{>0}$,)
$$\tR_-=-(ab-c)a'b(ab'c)\i\tS_+\in-\RR_{>0}$$
(since $\tS_+\in\RR_{>0}$).
We see that
$$\t\i(A,B)\sub\{(\tR_++a'/a,\tS_++b'/b),(\tR_-+a'/a,\tS_-+b'/b)\}$$
The same proof shows also the reverse inclusion so that
$$\t\i(A,B)=\{(\tR_++a'/a,\tS_++b'/b),(\tR_-+a'/a,\tS_-+b'/b)\}$$
Note that $(A,B)\m (\tR_++a'/a,\tS_++b'/b)$ is a continuous map
$\RR_{>0}^2@>>>\ct_{u,u'}$ which is the inverse of the continuous map
$\ct_{u,u'}@>>>\RR_{>0}^2$ defined by $\t$. This shows that 0.3 holds
in our case: one of the two maps $Z_i:U^-@>>>\RR$ ($i\in I$) attaches
to a matrix 
$$\left(\matrix 1&0&0\\x&1&0\\z&y&1\endmatrix\right)\in U^-$$
the number $z$; the other attaches to the same matrix the number $xy-z$.

\head 2. Proof of Proposition 0.5\endhead
\subhead 2.1\endsubhead
In this section we prove Proposition 0.5. 
For $u'\in U_{>0}^-,v\in U_{>0}^+$ we have $u'v\in G_{>0}$ hence we can
define $(u'v)^+\in U_{>0}^+$, $t_1\in T_{>0}$, $(u'v)^-\in U_{>0}^-$ by
the equation $u'v=(u'v)^+t_1(u'v)^-$. Note that the map
$U_{>0}^-\T U_{>0}^+@>>>U_{>0}^+\T T_{>0}\T U_{>0}^-$,
$(u',v)\m ((u'v)^+,t_1,(u'v)^-)$ is continuous.

\proclaim{Lemma 2.2} Let $B\in\cb_{>0}$. We write $B=u'B^+u'{}\i$ with
$u'\in U_{>0}^-$ uniquely determined. We have a homeomorphism
$$\{(v,t)\in U_{>0}^+\T T_{>0};t\in\ct_{(u'v)^-,u'}\}@>\si>>B\cap G_{>0}$$
given by $(v,t)\m u'vtu'{}\i$. 
\endproclaim
We have $B=\{u'vtu'{}\i;v\in U^+,t\in T\}$. If such $u'vtu'{}\i$ is in
$G_{>0}$ then $u'vt\in G_{>0}u'\sub G_{>0}U_{>0}^-\sub G_{>0}$ (see
\cite{L94, 2.12}) so that $u'vt=u'_0v_0t_0$ with
$u'_0\in U_{>0}^-,v_0\in U_{>0}^+,t_0\in T_{>0}$. It follows that
$u'=u'_0,v=v_0,t=t_0$ so that $v\in U_{>0}^+,t\in T_{>0}$.
Next we have $u'v=(u'v)^+t_1(u'v)^-$ as in 2.1) and
$$u'vtu'{}\i=(u'v)^+t_1(u'v)^-tu'{}\i
=(u'v)^+t_1t(t\i(u'v)^-tu'{}\i)\in U_{>0}^+T_{>0}U^-.$$
Since $u'vtu'{}\i\in G_{>0}$ we have also    
$u'vtu'{}\i=u_2t_2u'_2$ with $u_2\in U_{>0}^+,u'_2\in U_{>0}^-$,
$t_2\in T_{>0}$. It follows that $(u'v)^+=u_2,t_1t=t_2$,
$t\i(u'v)^-tu'{}\i=u'_2$. In particular, $t\i(u'v)^-tu'{}\i\in U_{>0}^-$.

Conversely, if $v\in U_{>0}^+$ and $t\in T_{>0}$ are given such that
$t\i(u'v)^-tu'{}\i\in U_{>0}^-$, then we have $u'vtu'{}\i\in G_{>0}$.
Indeed, as in 2.1 we have
$$\align&u'vtu'{}\i=(u'v)^+t_1(u'v)^-tu'{}\i\\&
=(u'v)^+(t_1t)(t\i(u'v)^-tu'{}\i)\in U_{>0}^+T_{>0}U_{>0}^-=G_{>0}.\endalign$$
The lemma follows.

\subhead 2.3\endsubhead
We prove 0.5. Let $B\in\cb_{>0}$. We write $B=u'B^+u'{}\i$ with
$u'\in U_{>0}^-$ uniquely determined. From 0.3 we obtain a homeomorphism 
$$\{(v,t)\in U_{>0}^+\T T_{>0};t\in\ct_{(u'v)^-,u'}\}@>>>
U_{>0}^+\T\RR_{>0}^I$$
given by $(v,t)\m(v,z((u'v)^-,u',t))$.
Composing this with the inverse of the homeomorphism in Lemma 2.2, we
obtain a homeomorphism of $B\cap G_{>0}$ with $U_{>0}^+\T\RR_{>0}^I$
so that 0.5 holds.

\head 3. Proof of Proposition 0.7 and of 0.8(c)\endhead
\subhead 3.1\endsubhead
Let $W$ be the Weyl group of $G$; let $\{s_i;i\in I\}$ be the simple
reflections in $W$. Let $w\m|w|$ be the length function; we have
$|s_i|=1$ for $i\in I$. For any $J\sub I$ let $W_J$ be the subgroup of
$W$ generated by $\{s_i;i\in J\}$ and let $w_0^J$ be the unique element
$w\in W_J$ with $|w|$ maximal. Let $w_0=w_0^I$. 

For $w\in W$ let $U^+(w)$ (resp. $U^-(w)$) be the subset of $U^+$ (resp.
$U^-$) defined in \cite{L94, 2.7} (resp. \cite{L94, 2.9}).

In the remainder of this section we fix $J\sub I$.
From the definitions we see that

(a) {\it $u\m uP_J^+u\i$ is a bijection $U^-(w_0w_0^J)@>\si>>\cp^J_{>0}$.}
\nl

\subhead 3.2\endsubhead
For any $\l=(\l_i)_{i\in I}\in\NN^I$ let $V_\l,\b_\l$ be as in 0.2.
Let $V_{\l,\RR}$ be the $\RR$-subspace of $V_\l$ spanned by $\b_\l$. 
Let $\PP$ be the set of all lines in $V_{\l,\RR}$. Note that $G$ acts naturally on $\PP$.
We shall assume that $\{i\in I;\l_i\ne0\}=I-J$. Then for
any $P\in\cp^J$ there is a unique $L_P^\l\in\PP$ such that the stabilizer
of $L_P^\l$ in $G$ is equal to $P$; moreover, 

(a) {\it the map $P\m L_P^\l$ from $\cp^J$ to $\PP$ is injective.}
\nl
A line in $\PP$ is said to be in $\PP_{>0}$ if it is spanned by a linear combination of
elements in $\b_\l$ with all coefficients being in $\RR_{>0}$.
We shall now assume in addition that $\l$ is such that
$\l_i$ is sufficiently large for any $i\in I-J$ so that \cite{L98, 3.4} is applicable;
thus the following holds:

(b) {\it For $P\in\cp^J$ we have $P\in\cp^J_{>0}$ if and only if $L_P^\l\in\PP_{>0}$.}
\nl
Now let $g\in G_{>0}$. From \cite{L94, 5.2} we see that the following holds.

(c) {\it There is a unique line $L_g\in\PP_{>0}$ such that $gL_g=L_g$.}
\nl

\subhead 3.3\endsubhead  
Let $g\in G_{>0}$. We prove Proposition 0.7.

By \cite{L94, 8.9} there exists $B\in\cb_{>0}$ such that $g\in B$. Let $P\in\cp^J$
be such that $B\sub P$. We have $P\in\cp^J_{>0}$, $g\in P$. This proves the existence
in 0.7. Assume now that $P'\in\cp^J_{>0},P''\in\cp^J_{>0}$ satisfy
$g\in P',g\in P''$. In the setup of 3.2(b) we have 
 $L_{P'}^\l\in\PP_{>0}$,  $L_{P''}^\l\in\PP_{>0}$. Since $g\in P'$, the line
 $L_{P'}^\l$ is $g$-stable hence, with notation of 3.2(c), we have $L_{P'}^\l=L_g$.
Similarly we have $L_{P''}^\l=L_g$. Thus we have $L_{P'}^\l=L_{P''}^\l$. Using 3.2(a)
we deduce $P'=P''$. This proves the uniqueness in 0.7.

\subhead 3.4\endsubhead
We now fix $P\in\cp^J_{>0}$. By 3.1(a) we have $P=uP_J^+u\i$ for a well defined
$u\in U^-(w_0w_0^J)$.

We show:

(a) {\it $v\m uvB^+v\i u\i$ is a bijection $U^-(w_0^J)@>\si>>\g_P$. }
\nl
If $v\in U^-(w_0^J)$ then $uv\in U^-(w_0)=U^-_{>0}$ hence 
$uvB^+v\i u\i\in\cb_{>0}$; moreover we have $vB^+v\i\in P_J^+$ hence
$uvB^+v\i u\i\sub uP_J^+u\i=P$ so that $uvB^+v\i u\i\in\g_P$. Thus (a) is a well defined
map. Now let $B\in\g_P$. We have $B=u_1B^+u_1\i$ where $u_1\in U^-(w_0)$ and
$B\sub uP_J^+u\i$ that is $u_1B^+u_1\i\sub uP_J^+u\i$. Now there is a unique $P'\in\cp^J$
containing $u_1B^+u_1\i$. Since $u_1P_J^+u_1\i$ and $uP_J^+u\i$ are such $P'$, we must have
$u_1P_J^+u_1\i=uP_J^+u\i$.
We have $u_1=u'_1u''_1$ where $u'_1\in U^-(w_0w_0^J)$, $u''_1\in U^-(w_0^J)$;
Hence $u_1P_J^+u_1\i=u'_1u''_1P_J^+u''_1{}\i u'_1{}\i=u'_1P_J^+u'_1{}\i$ so that
$u'_1P_J^+u'_1{}\i=uP_J^+u\i$. Using this and 3.1(a) we deduce $u'_1=u$ so that 
$B=uu''_1B^+u''_1{}\i u\i$. We see that the map (a) is surjective.
The injectivity of the map (a) is immediate. This proves (a).

From the definitions we see that (a) is a homeomorphism. Hence 0.8(c) holds.

\subhead 3.5\endsubhead
Let $P\in\cp^J_{>0}$ and let $P_{red}$ be the reductive quotient of $P$.
Define $u\in U^-(w_0w_0^J)$ by $P=uP_J^+u\i$ (see 3.1(a)).
Now $P_{red}$ has a natural pinning. (It is induced by the obvious pinning of the
reductive quotient of $P_J^+$ using the isomorphism $P_J^+@>>>P$ given by conjugation by 
$u$.) Hence we can define $\cb_{P,>0}$ (an open subset of the real flag manifold of 
$P_{red}$) in terms of $P_{red}$ with its pinning in the same way as $\cb_{>0}$ was defined
in terms of $G$ with its pinning. (Although $P_{red}$ is not necessarily semisimple, the
same definitions can be applied.) We have a bijection

(a) $\cb_{P,>0}@>\si>>\g_P$
\nl
obtained by taking inverse image under the obvious map $P@>>>P_{red}$.
This can be deduced from the proof in 3.4.

\enddocument